\documentclass[a4paper,12pt]{article}
\usepackage{amssymb}

\title{$L^p$--distributions on symmetric spaces
}
\author{Michael Ruzhansky}
\date{}
\begin{document}

\begin{flushleft}
appeared in: Results Math.  44 (2003), 159--168.
\end{flushleft}

\vspace{0.7cm}

\begin{center} 
{\LARGE $L^p$--distributions on symmetric spaces}
\end{center}
\medskip

\begin{center} 
{\large Michael Ruzhansky}
\end{center}
\medskip


 \newtheorem{prop}{Proposition}
 \newtheorem{lemma}{Lemma}
 \newtheorem{cor}{Corollary}
 \newtheorem{theorem}{Theorem}
 \newtheorem{definition}{Definition}
 \newtheorem{remark}{Remark}

  \newcommand{\Rn}{{{\mathbb R}^n}}
  \newcommand{\Rone}{{\mathbb R}}
  \newcommand{\Zed}{{\mathbb Z}}
  \newcommand{\Zpos}{{\mathbb Z}_{\geq 0}}
  \newcommand{\Nat}{{\mathbb N}}
  \newcommand{\Dmsy}{{\mathbb D}}
  \newcommand{\Cmplx}{{\mathbb C}}
  \newcommand{\DLprime}{{{\cal D}_L^\prime}}
  \newcommand{\Aloc}{{\cal A}_{loc}}
  \newcommand{\Dcal}{{\cal D}}
  \newcommand{\Dlp}{{{\cal D}_{L^p}}}
  \newcommand{\Dlr}{{{\cal D}_{L^r}}}
  \newcommand{\Dli}{{{\cal D}_{L^\infty}}}
  \newcommand{\Dlo}{{{\cal D}_{L^1}}}
  \newcommand{\Dlid}{{{\dot{\cal D}}_{L^\infty}}}
  \newcommand{\Dlq}{{{\cal D}_{L^q}}}
  \newcommand{\Dcalpr}{{{\cal D}^\prime}}
  \newcommand{\Dlopr}{{{\cal D}^\prime_{L^1}}}
  \newcommand{\Dlppr}{{{\cal D}^\prime_{L^p}}}
  \newcommand{\Dlrpr}{{{\cal D}^\prime_{L^r}}}
  \newcommand{\Dlqpr}{{{\cal D}^\prime_{L^q}}}
  \newcommand{\Acal}{{\cal A}}
  \newcommand{\Bcal}{{\cal B}}
  \newcommand{\Scal}{{\cal S}}
  \newcommand{\Ecal}{{\cal E}}
  \newcommand{\Ecalm}{{{\cal E}^{(m)}}}
  \newcommand{\Fcal}{{\cal F}}
  \newcommand{\Ecalpr}{{\cal E}^\prime}
  \newcommand{\Ccal}{{\cal C}}
  \newcommand{\contcc}{{\cal C}_c}
  \newcommand{\Bdot}{\dot{{\cal B}}}
  \newcommand{\Mb}{{\cal M}_b(\Rn)}

  \newcommand{\HPoincare}{{\mathbb H}}
  \newcommand{\ddx}{\frac{\partial}{\partial x}}
  \newcommand{\ddy}{\frac{\partial}{\partial y}}
  \newcommand{\dfdx}{\frac{\partial f}{\partial x}}
  \newcommand{\dfdy}{\frac{\partial f}{\partial y}}
  \newcommand{\ddxc}[1]{\frac{\partial #1}{\partial x}}
  \newcommand{\ddyc}[1]{\frac{\partial #1}{\partial y}}

  \newcommand{\ddxj}[1]{\frac{\partial #1}{\partial x_j}}

  \newcommand{\Ecpr}{{\cal E}^\prime(\Rn)}
  \newcommand{\Dcpr}{{\cal D}^\prime(\Rn)}


\begin{abstract}
    The notion of $L^p$--distributions is introduced on
    Riemannian symmetric spaces of noncompact type and
    their main properties are established.  We use a
    geometric
    description for the topology of the space of test functions
    in terms of the Laplace--Beltrami operator. 
    The techniques are based
    on a-priori estimates for elliptic operators. We show
    that structure theorems, similar to $\Rn$, hold on
    symmetric spaces. We give estimates for the convolutions.
    
\end{abstract}

%
%

 \section{Introduction}

 In this paper we will generalize the notion of $L^p$--distributions
 to the setting of symmetric spaces and establish their
 main properties. Basic examples are $L^p$-functions and their generalized
 derivatives.
 We will give an invariant description of these distributions
 as sums of iterated Laplace-Beltrami operators applied to 
 $L^p$-functions.  As consequence, we give several results on 
 convolution and other properties
 of these distributions.
 
 We use a geometric approach
 based on a-priori estimates involving elliptic operators.
 We suggest alternative
 definitions of test functions using the iterated powers of
 the Laplace--Beltrami operator associated to the Riemannian structure.
 Combining the regularity properties of pseudo-differential
 operators with the existence of fundamental solutions
 for elliptic invariant operators on a symmetric space,
 we can obtain certain uniform a-priori estimates for
 the elliptic invariant operators.
 These estimates are established in
 $L^p$ spaces for $1\leq p \leq\infty$.
 Similar estimates for the gradient
 are known for $p=\infty$ on manifolds
 with bounded curvature (\cite{JK}). On the other hand, estimates
 for $1<p<\infty$ are well known for Euclidean spaces
 (\cite{Stein}) and are
 closely related to the continuity
 of pseudo-differential operators of order zero, which holds only
 locally on general manifolds.
 We will give a unifying algebraic proof of global
 a-priori estimates in
 $L^p$--space for all $1\leq p\leq\infty$ on symmetric
 spaces of the noncompact type. The estimates are uniform
 since one can use convolution to extend them globally.
 The general theory of the second order differential
 operators on Lie groups can be found in \cite{Rob}, \cite{Var},
 and we refer there for more detailed information.

 In the last section we use these estimates to define the spaces
 of test functions dependent on $p$. Then the standard construction
 leads to the distribution spaces $\Dlppr$.
 We relate these spaces to the
 corresponding $L^p$ spaces and their generalized derivatives.
 Distribution spaces of this type are useful in a number of
 applications (\cite{Schwartz}, \cite{BN}). The general distribution
 theory on symmetric spaces can be found in \cite{Helg1} 
 (as well as 
 many general results on invariant
%
%
  differential operators
 and convolutions).
 However, in comparison with $\Rn$, we obtain the representation
 formulas involving the iterated Laplacian.
 Finally we show that the convolution properties of
 \cite{Schwartz} hold.
 The spaces $\Dlppr$ provide a scale of distribution spaces
 which leads to the tempered distributions. In \cite{Cator},
 summable distributions ($\Dlopr$) are used as a foundation for the
 theory of distributions in locally convex spaces. In the
 case of symmetric spaces, this leads to the distribution
 theory as well, since the convolution techniques are available.
 We hope to develop this point of view in subsequent work.

 Finally we note that in this paper we consider groups with classical
 Lie algebras only. The exceptional cases
 ${\mathfrak e}_6$, ${\mathfrak e}_7$, ${\mathfrak e}_8$ are left
 out because the mapping $d\pi$ induced by the canonical projection
 $\pi$ of the symmetric space is not surjective from the center
 of the algebra of left invariant differential operators
 on the group to the algebra of left invariant differential
 operators on the symmetric space (\cite{Helg3}).

 I would like to thank professor Erik Thomas for drawing my attention to
 the subject and for clarifying the equicontinuity argument
 in Theorem \ref{th:thconv}.

\section{Preliminary estimates}

 Let $M=G/H$ be a Riemannian symmetric space of the noncompact type.
 This means, that it can be viewed as a quotient  $M=G/H$, where
 $G$ is a connected semisimple Lie group with trivial center and
 $H$ is its maximal compact subgroup.
 We will always assume that the Lie algebra ${\mathfrak g}$ of
 $G$ is classical. Therefore, we exclude the exceptional cases
 ${\mathfrak e}_6$, ${\mathfrak e}_7$, ${\mathfrak e}_8$, with
 pairs $(G,H)$ listed in \cite{Helg3}.
 The Riemannian structure on $M$
 is supposed to be invariant under the left action of $G$.

We will
describe first a relation between parametrix and fundamental
solutions of $P$. This
will lead to some integrability properties of fundamental solutions,
which will be applied for getting uniform $L^p$ estimates.
Let $P$ be an elliptic differential operator of order $m$ on $M$.
Then there exists a parametrix for $P$, namely a pseudo-differential
operator $Q\in\Psi^{-m}(M)$ of order $-m$, such that
  \begin{equation}
    QP=I+R
    \label{eq:parman}
  \end{equation}
 with
 $R\in\Psi^{-\infty}(M)$.
Let $F$ be the distributional kernel of $Q$. This means that
$F\in \Dcalpr(M\times M)$ satisfies
$\langle Q\phi,\psi\rangle =\langle F,\phi\otimes\psi\rangle$
for all $\phi,\psi\in\Dcal(M)$, where $\Dcal(M)$ is the standard
space of smooth compactly supported functions on $M$.
In the sequel we will write

  $$ (Qu)(x) = \int_M F(x,y) u(y) d\mu(y),$$
and we view $Q$ as a singular integral operator with kernel $F$,
where $\mu$ is the Riemannian measure on $M$.
Note also, that $P$ is
automatically proper supported because it is a local operator.
Assume now that $P$ is $G$--left invariant. Then, according to
\cite[Theorem 4.2]{Helg2}, there exists a fundamental solution for $P$, namely
a distribution $K\in\Dcalpr(M)$, such that

\begin{equation}
      PK=\delta,
   \label{eq:fs1}
\end{equation}
$\delta$ being a delta function at the origin $p$ of $M$.
Let $\alpha\in\Dcal(M)$ be a test function, such that
$\alpha(x)=1$ for $x$ in a small neighborhood of $p$.
Then equality (\ref{eq:fs1}) implies the existence of $\beta\in\Dcal(M)$,
such that

\begin{equation}
      P(\alpha K)=\delta+\beta.
   \label{eq:fs2}
\end{equation}
In fact, taking $\beta=P(\alpha K)-\delta$ one readily verifies
that $\beta\in\Dcal(M)$.
The only singularity of $K$ occurs at the point $p$ due to the ellipticity
of $P$ and we will be interested in the integrability properties of
$\alpha K$.
An application of formula (\ref{eq:parman}) to $\alpha K$ yields

\begin{equation}
  \alpha K + R(\alpha K) = QP(\alpha K)  = Q\delta  + Q\beta,
  \label{eq:rel1}
\end{equation}
the last equality due to (\ref{eq:fs2}).
We have $\beta\in\Dcal(M)$ and $R\in\Psi^{-\infty}(M)$ implying
$R(\alpha K), Q\beta\in\Ccal^\infty(M)$. The operators $Q$ and $R$ are
properly
supported, therefore all the functions in (\ref{eq:rel1}) have compact
support. Let $D\in\Psi^k(M)$ be properly supported, $k<m$.
The application of $D$ to
(\ref{eq:rel1}) and the arguments above imply

\begin{equation}
  D(\alpha K) =  DQ\delta + \psi,
  \label{eq:rel2}
\end{equation}
with $\psi\in\Dcal(M)$. Now, the operator $DQ$ is of a negative order
$k-m$ and
this implies the integrability of of its integral
kernel at $p$. This property holds locally on arbitrary smooth manifolds
(cf. \cite{Stein}, and some related properties can be found in
\cite{RuzhS}).
The latter is equal to $DQ\delta(x)$.
Equality (\ref{eq:rel2}) implies the integrability of $D(\alpha K)$.
Thus, we have proved the following

\begin{lemma}
   Let $P$ be an invariant elliptic differential operator of order
   $m$ on $M$
   and $K$ its fundamental solution at $p$. Then, for every
   $D\in\Psi^k(M)$, $k<m$, $DK$ is locally integrable
   $DK\in L^1_{loc}(M)$.
  \label{l:lfs}
\end{lemma}
We will apply this lemma to two cases, $D$ being an invariant
differential operator and $P$ being the Laplace--Beltrami operator on
$M$ equipped with a Riemannian structure.
The space of $G$--left invariant differential operators of order $k$
on $M$ will be denoted $\Dmsy^k(M)$.
Let $\Zed(G)$ denote the center of the algebra of the left invariant
differential operators on $G$. Let $\pi:G\to M=G/H$ be the canonical
projection and let ${\mathfrak g}$ denote the Lie algebra of $G$.
Note, that
$d\pi : {\mathfrak g}\to T_pM$ can be extended to the algebra
$\Dmsy(G)$ of the left
invariant differential operators on $G$.

\begin{theorem}
  Let the Lie algebra ${\mathfrak g}$ of $G$ be classical and semisimple.
  Let $P\in\Dmsy^m(M)$ be elliptic, $D\in\Dmsy^k(M)$, $0<k<m$ and
  $1\leq p\leq\infty$.
  Then there exist constants $A,\; B$, such that for every
  $u\in L^p(M)$ with $Pu\in L^p(M)$, we have $Du\in L^p(M)$
  and
    \begin{equation}
       ||Du||_p \leq A ||Pu||_p + B ||u||_p.
      \label{eq:estsym1}
    \end{equation}
  If $p=\infty$, then $Du$ is continuous.
  \label{th:th1}
\end{theorem}
If ${\mathfrak g}$ is not classical (for example, when ${\mathfrak g}$
is the real form of the exceptional Lie algebras
${\mathfrak e}_6$, ${\mathfrak e}_7$, ${\mathfrak e}_8$),
the statement of Theorem \ref{th:th1} still holds if we assume
that $P\in d\pi(\Zed(G))$.

 Let $X$ be a smooth
 vector field on $M$. $X$ is called {\em bounded} if
 there exist a constant $C$ such that  $||X_x||\leq C$
 for every point
 $x\in M$, where $||\cdot ||=\langle\cdot,\cdot\rangle^{1/2}$
 is the Riemannian norm on $T_xM$, corresponding to the Riemannian structure.

\begin{theorem}
  Let $M$ be a Riemannian symmetric space as above
  and let $\Delta$ be the associated
  Laplace--Beltrami operator. Let $1\leq p\leq\infty$ and let $X$ be
  a smooth bounded
  vector field on $M$.
  Then there exist constants $A,\; B$, such that for every
  $u\in L^p(M)$ with $\Delta u\in L^p(M)$, we have that the derivative
  of $u$ with respect to $X$ is $L^p$--integrable, $Xu\in L^p(M)$
  and
    \begin{equation}
       ||Xu||_p \leq A ||\Delta u||_p + B ||u||_p.
      \label{eq:estsym2}
    \end{equation}
  If $p=\infty$, then $Xu$ is continuous. Moreover, $A$ and $B$ can be
  chosen independently over the set of the smooth vector fields $X$
  bounded by 1.
  \label{th:th2}
\end{theorem}
  The statements of Theorems \ref{th:th1} and \ref{th:th2} can be
  improved with respect to the order $k$. However, we will not
  discuss it here since we do not need it for the applications in
  the next section.

To show the statements of the theorems,
first, convolving (\ref{eq:fs2}) with $u\in\Dcalpr(M)$, we conclude
that

\begin{equation}
  u = u*\delta =  u* P(\alpha K) - u*\beta.
   \label{eq:idrepr1}
\end{equation}
For classical algebras ${\mathfrak g}$ the extension
$d\pi : {\mathfrak g}\to T_pM$ to the algebra of the left
invariant differential operators on $G$ has the property
$d\pi(\Zed(G))=\Dmsy(M)$, (\cite[Proposition 7.4, Theorem 7.5]{Helg2}
or \cite[Remark, p.326]{Helg1}).
This means that
$P\in\Dmsy^m(M)$ is an image of a bi-invariant operator on $G$ and,
in particular, commutes with left and right convolution.
Therefore, (\ref{eq:idrepr1}) implies

\begin{equation}
   u = Pu * \alpha K - u * \beta.
   \label{eq:idrepr2}
\end{equation}
Let us make several remarks.
In general, for a separable unimodular Lie group $G$ and
a compact $H$, such that $(G,H)$ is a symmetric pair,
equation (\ref{eq:idrepr2}) follows from (\ref{eq:idrepr1}), for
an arbitrary $P\in\Dmsy(G/H)$. Indeed, distribution $\alpha K$
is compactly supported and, according to
\cite[Theorem 5.5, p. 293]{Helg1}, we have
$u*P(\alpha K)=Pu * \alpha K.$
In our case, (\ref{eq:idrepr2}) follows from (\ref{eq:idrepr1})
in conditions of Theorem \ref{th:th1}
because $P$ is the image of a bi-invariant operator on $G$.
Note, that if the Lie group $G$ is compact, then statements of
Theorem \ref{th:th1}
and Theorem \ref{th:th2}
follow directly from the local estimates of the same
type. Indeed, if we apply equality (\ref{eq:parman}) to $u$ and then
apply $D$ (not necessarily invariant) to it, we get
$$Du=DQPu-DRu.$$
Now, operators $DQ$ and $DR$ are pseudo-differential operators of strictly
negative orders and, therefore, locally bounded on $L^p$.
Therefore, we may assume that $G$ is not compact. Applications of
the next section will deal with the case of the iterated
Laplace--Beltrami operator $P=\Delta^l$. According to
\cite[p.331]{Helg1}, if the Riemannian structure on $G$ is
associated to the Killing form, the Laplace--Beltrami operator on $M=G/H$
is the image of the Casimir operator on $G$, which belongs to
the center $\Zed(G)$.

\bigskip
\noindent
{\bf Proof of Theorem \ref{th:th1}:}
An application of $D$ to equality (\ref{eq:idrepr2}) together with
an argument above imply
$$ Du = Pu * D(\alpha K) - u * D\beta. $$
By Lemma \ref{l:lfs} we have $D(\alpha K)\in L^1(M)$ and Young
inequality (\cite{Hew})
yields estimate (\ref{eq:estsym1}). In case $p=\infty$
we have convolutions of the type $L^\infty * L^1$, which give
continuous functions. This completes the proof of Theorem
\ref{th:th1}. \\


Now we will need some notation and auxiliary results.
We start with constructing an invariant Riemannian structure
on $G$, making the projection $\pi$ a Riemannian submersion.
In the sequel $l_g h = gh, \; r_g h = hg, \; g,h\in G$
will denote the left and right group actions on $G$. The induced
actions of $G$ on $M$ will be denoted by the same letters.
The adjoint
representation will be denoted by ${\rm Ad}:G\to GL(T_eG)$.
The following lemma is standard.

\begin{lemma}
There exists a $G$--left and $H$--right invariant Riemannian structure
on $G$, such that the canonical projection $\pi:G\to M$ is a Riemannian
 submersion, i.e. a submersion, for which
 horizontal lift of vector fields preserves  Riemannian norms.
 \label{l:Riem}
\end{lemma}
%
 For every $g\in G$ let $K_g={\rm Ker}\; d_g\pi$
 ($K_g \cong T_gH \cong T_eH$).
 Let $N_g$ be the orthogonal complement to $K_g$ with respect to
 ${\rm Ad}(H)$-invariant inner product
 $(\cdot,\cdot)_0$ given by
  \begin{equation}
   (X,Y)_0 := \int_{{\rm Ad}(H)} \langle A(X), A(Y) \rangle_0 d\mu(A),
   \label{eq:inv}
  \end{equation}
 where $\mu$ is Haar measure on a compact
 set ${\rm Ad}(H)$ and $\langle\cdot,\cdot\rangle_0$ is an
 arbitrary inner product on $T_eG$.
 Thus we have $T_g G=K_g \bigoplus N_g$.
 Define an inner product $(\cdot,\cdot)_{N_e}$
 on $N_e$ for vectors $\bar{X},\bar{Y}\in N_e$ by

  \begin{equation}
    (\bar{X},\bar{Y})_{N_e} := \langle X,Y\rangle_{M_p},
    \label{eq:invn}
  \end{equation}
 where $\langle\cdot,\cdot\rangle_{M_p}$ is the restriction to $T_pM$ of
 the given invariant Riemannian structure on $M$
 and $X=d\pi(\bar{X}),Y=d\pi(\bar{Y})\in T_{p}M$.
 Vectors $\bar{X}$ and $\bar{Y}$ are uniquely defined, $d\pi$ being
 an isomorphism between $N_e$ and $T_pM$, and they are
 called the {\em horizontal lifts} of
 $X$ and $Y$.
 The desired Riemannian structure on $G$ can now be constructed
 by applying $dl_g$ to
   \begin{equation}
    \langle u, v \rangle :=
      (u|_{K_e},v|_{K_e})_0 + (u|_{N_e},v|_{N_e})_{N_e},
   \label{eq:invghm}
  \end{equation}
 where $u,v\in T_eG$, $u|_{K_e},v|_{K_e}$ and $u|_{N_e},v|_{N_e}$
  are projections
 of $u,v$ on $K_e$ and $N_e$, respectively.
 The inner product in
 (\ref{eq:invghm}) is clearly ${\rm Ad}(H)$--invariant,
 the expansion is therefore $G$--left and $H$--right invariant.
 It follows immediately from formulas
 (\ref{eq:invn}) and (\ref{eq:invghm}) that all
 $d_g\pi$ are partial isometries (isometries from $N_g$ to
 $T_{\pi(g)}M$). We fix this inner product on $G$ since it
 satisfies Lemma \ref{l:Riem}.
 For the sake of completeness let us list briefly some well known
 properties
 of the pullback $\sharp$ which will be necessary.

 \begin{lemma}
   \begin{itemize}
   \item[{\rm (i)}]
      Let $\phi\in\contcc(G)$ be a continuous
  compactly supported function on $G$.
  Then for $x=\pi(g)$ the function $\phi^\flat(x)$ is correctly
  defined by
    $$ \phi^\flat(x)=\int_H \phi(gh)dh,$$
     where
  $dh$ is the normalized Haar measure on $H$.
  Moreover,  $\phi^\flat\in \contcc(M)$ and
  mapping $\phi\to\phi^\flat$ is linear surjective from $\contcc(G)$
  to $\contcc(M)$ and from $\Dcal(G)$ to $\Dcal(M)$.

 \item[{\rm (ii)}]
     The transpose of $\phi\to\phi^\flat$ defined by
   $$\langle T^\sharp,\phi \rangle = \langle T,\phi^\flat \rangle$$
     is an injective mapping from $\Dcalpr(M)\to\Dcalpr(G)$.

 \item[{\rm (iii)}]
       Let $S\in\Dcalpr(G)$. Then
     $S$ is right $H$-invariant if and only if there exists
     $T\in\Dcalpr(M)$, such that $S=T^\sharp$.
     For $T_1, T_2\in\Dcalpr(M)$ the convolution products on
     $G$ and $M$ are related by
     $$ T_1^\sharp * T_2^\sharp = (T_1 * T_2)^\sharp.$$

 \item[{\rm (iv)}]
          Let $Y$ be a horizontal lift of a vector field $X$ on $M$
    and let $T\in\Dcalpr(M)$ be a distribution on $M$. Then
    $Y(T^\sharp)=(XT)^\sharp$.
 \end{itemize}
   \label{l:facts}
 \end{lemma}
{\bf Proof of Theorem \ref{th:th2}:}
 The pullback of formula (\ref{eq:idrepr2}) now reads

 \begin{equation}
   u^\sharp =
  (\Delta u)^\sharp * (\alpha K)^\sharp - u^\sharp * \beta^\sharp.
   \label{eq:gidrepr}
 \end{equation}
 Let $X$ be a smooth vector field on $M$, bounded by one:
 $||X_x||\leq 1$.
  Let $Y$ denote the horizontal lift of $X$:
  \begin{itemize}
    \item[1.] $d_g\pi(Y_g)=X_x$, where $x=\pi(g)$.
    \item[2.] $Y_g\in N_g=({\rm Ker}\; d_g\pi)^\perp$.
  \end{itemize}
 It is not difficult to see that $Y$ is smooth.
   Let $Y_1,\ldots,Y_N$ be an orthonormal basis
   of the Lie algebra {$\mathfrak g$},
  such that $(Y_1)_g,\ldots,(Y_n)_g\in N_g$ for all $g\in G$.
  Vector field $Y$ can be decomposed with respect to the basis
  $Y_1,\ldots,Y_n$ at every point $g\in G$:

    \begin{equation}
       Y_g = \sum_{i=1}^n a_i(g) Y_{i,g}\in N_g\subset T_gG,
       \label{eq:decyyi}
    \end{equation}
  where $Y_{i,g}=(Y_i)_g=d_el_g (Y_i)_e$ are values at $g$
  of the left invariant vector fields $Y_i$.
  Note, that such decomposition is only pointwise because $Y$ need not be
  left invariant in general, we use that $Y_g\in N_g$ and the fact that
  $Y_{1,g},\ldots,Y_{n,g}$ constitute a basis for a linear space $N_g$.
  However, it is global and
  functions $a_1,\ldots,a_n$ are smooth due to the smoothness of $Y$
  and $Y_1,\ldots,Y_n$.

  The norm of $Y_g$ at $T_gG$ is
  $||Y_g||^2=\sum_{i=1}^n |a_i(g)|^2$.  In view
  of Lemma \ref{l:Riem}, $||Y_g||=||X_x||\leq 1$.
  In particular, $|a_i(g)|\leq 1$ for all $g\in G$.
  Now we differentiate $u^\sharp$ in
  (\ref{eq:gidrepr}) with respect to the basis
  vector fields $Y_i$ and the left
  invariance of $Y_i$ yields:

    \begin{equation}
      Y_i u^\sharp = (\Delta u)^\sharp * Y_i(\alpha K)^\sharp +
      u^\sharp * Y_i \beta^\sharp.
      \label{eq:yi1}
    \end{equation}
  Obviously $Y_i \beta^\sharp \in \Dcal(G) \subset L^1(G)$.
  In view of Lemma \ref{l:lfs} the compactly supported distribution
  $\alpha K$ and its derivatives are integrable and so are their
  pullbacks, the pullback mapping being an isometry of $L^p$ spaces.
  Let $A_i = ||Y_i(\alpha K)^\sharp||_1$ and
      $B_i = ||Y_i \beta^\sharp||_1$.
  Application of Young inequality \cite[Cor. 20.14]{Hew}
  to (\ref{eq:yi1}) yields
      $$ ||Y_i u^\sharp||_p \leq A_i ||(\Delta u)^\sharp||_p +
                               B_i ||u^\sharp||_p. $$
  Decomposition (\ref{eq:decyyi}) together with bounds on $a_i$
  and equalities   $||(\Delta u)^\sharp||_p = ||\Delta u||_p$ and
  $||u^\sharp||_p = ||u||_p$ imply
           $$ ||Y u^\sharp||_p \leq A ||\Delta u||_p +
                                   B ||u||_p$$
  with $A = \sum_{i=1}^n A_i$ and  $B = \sum_{i=1}^n B_i$.
  By Lemma \ref{l:facts} we have $||Yu^\sharp||_p=||Xu||_p$,
  establishing inequality (\ref{eq:estsym2}).
  In case $p=\infty$, formulas (\ref{eq:yi1}) and (\ref{eq:decyyi})
  imply the continuity of $Yu^\sharp$. By Lemma \ref{l:facts},
  $(Xu)^\sharp=Yu^\sharp$ is continuous. The continuity of $Xu$
  now follows from the fact that $M$ is equipped with the quotient
  topology, i.e. the strongest topology, for which $\pi$ is a
  continuous mapping. This finishes the proof of Theorem \ref{th:th2}.

 \section{$L^p$--distributions}

 Let $M$ be a symmetric space as before. In this section we will apply
 estimates of the previous section to the particular cases of $P=\Delta^k$,
 to construct spaces $\Dlppr (M)$ of $L^p$--distributions, invariantly
 on $M$.
 For $1\leq p\leq\infty$, we consider the space
$$ \Dlp (M) =\{ \phi\in C^\infty(M): \; \Delta^k \phi\in L^p(M),\;
        \;\forall\; k\in\Zpos\},$$
equipped with a countable system of seminorms
$$ \omega_{p,k}(\phi)=\max\{||\phi||_p,||\Delta^k\phi||_p\},\; k\in\Zpos.$$
They define the coarsest locally convex topology for which the
maps $\Delta^k:\Dlp(M)\to L^p(M)$ are continuous for all $k\in\Zpos$.
With this topology $\Dlp (M)$ become Fr\'echet spaces.
For $1\leq p<\infty$ we obviously have the continuous embeddings
$$ \Dcal(M)\subset \Dlp(M)\subset \Dcalpr(M).$$
The first inclusion is also dense for $p<\infty$. For $p=\infty$
define $\Dlid$ to be the subspace of $\Dli$ of functions vanishing
at infinity. Then $\Dcal(M)$ is dense in $\Dlid$, so the spaces
$\Dlp$ for $1\leq p<\infty$ and $\Dlid$ are normal spaces of distributions.
This implies that their strong duals are the subspaces of $\Dcalpr (M)$.
For $1<p\leq\infty$ and $q$ such that $1/p+1/q=1$, we denote by
$\Dlppr (M)$ the strong dual of $\Dlq (M)$. For $p=1$ we denote
by $\Dlo^\prime(M)$ the strong dual of $\Dlid (M)$. In view of Theorem
\ref{th:th2} the spaces $\Dlp (M)$ coincide with the spaces $\Dlp (\Rn)$
of L. Schwartz (\cite[Ch. VI]{Schwartz}) when $M=\Rn$.

\begin{theorem} A distribution $T$ belongs to $\Dlppr (M)$ if
and only if there exist $m=m(T)\in\Nat$ and $f,g\in L^p(M)$,
such that
   \begin{equation}
       T= f+\Delta^m g.
     \label{eq:str1}
   \end{equation}
  \label{th:threpr1}
\end{theorem}
{\bf Proof:} Let $1/p+1/q=1$.
First, if $T$ is given by formula (\ref{eq:str1}), then $T$ is
a linear continuous functional on $\Dlq (M)$, so $T\in\Dlppr(M)$.

Conversely, suppose that $T\in\Dlppr(M)$ and assume that $p>1$.
Then there exists a number $m$ such that
  \begin{equation}
   |\langle T,\phi\rangle|\leq C\omega_{q,m}(\phi),
   \label{eq:str1eq1}
  \end{equation}
for every $\phi\in\Dlq(M)$.
Let $\iota:\Dlq (M)\to L^q(M)\times L^q(M)$ be an injective inclusion,
$\iota(\phi)=(\phi,\Delta^m\phi)$. On the image of $\iota$ consider
a linear map $L:\iota(\Dlq(M))\to\Cmplx$ defined by
  $$ L(\phi,\Delta^m\phi)=\langle T,\phi\rangle.$$
Inequality (\ref{eq:str1eq1}) implies
  $$ |L(\phi,\Delta^m\phi)|\leq C\omega_{q,m}(\phi)=
        C\max\{||\phi||_q,||\Delta^m\phi||_q\},$$
which means that $L$ is continuous if we equip $\iota(\Dlq(M))$ with
the induced topology of $L^q(M)\times L^q(M)$. Therefore, by Hahn--Banach
theorem $L$ allows a linear continuous extension to $L^q(M)\times L^q(M)$,
which we also denote by $L$. Now, the dual
of $L^q(M)\times L^q(M)$ is $L^p(M)\times L^p(M)$, implying
 $$\langle T,\phi \rangle=L(\phi,\Delta^m\phi)=
    \int_M\phi f d\mu+\int_M (\Delta^m\phi)g\, d\mu =$$
 $$ \langle f,\phi\rangle + \langle g, \Delta^m\phi\rangle=
    \langle f,\phi\rangle + \langle \Delta^m g,\phi\rangle=$$
 $$ \langle f+\Delta^m g,\phi \rangle$$
 for all $\phi\in \Dlq(M)$, and some $f,g\in L^p(M)$, with integration
 with respect to the Riemannian measure $\mu$.
The case with $p=1$ follows the same lines, but one should take
$C_0(M)\times C_0(M)$ instead of $L^q(M)\times L^q(M)$.
Then $f$ and $g$ are Radon measures, but
a standard closed graph argument shows that $f$ and $g$ can be
taken in $L^\infty$.
The proof is complete.

From Theorem \ref{th:threpr1} we can imagine the general structure
of $L^p$ distributions. For example, if $f\in L^p(M)$ is compactly
supported and $g\in L^p(M)$ is supported ``near infinity", then
sums of derivatives of $f$ and $g$ represent elements of $\Dlppr(M)$.
There is an improvement to Theorem \ref{th:threpr1}:

\begin{cor} Let $T\in\Dlppr(M)$. Then there exist $m=m(T)\in\Nat$
 and continuous functions $f_k\in L^p(M)$, $0\leq k\leq m$, such that
  $$ T=\sum_{k=0}^m \Delta^k f_k.$$
\end{cor}
{\bf Proof:} Let $Q_j$ be a right parametrix for the elliptic operator
  $(I-\Delta)^j$. Then we have $(I-\Delta)^j Q_j=I+R_j$, with $R_j$
 a smoothing operator. Applying this to a function $h\in L^p$, we
 get
  $$ h=(I-\Delta)^j (Q_j h)-R_j h.$$
 Choosing $j$ sufficiently large, the functions $Q_j h$ and $R_j h$
 are continuous $L^p$-functions.  This argument for functions
 $f$ and $g$ in (\ref{eq:str1}) imply the corollary.

\begin{theorem}
   A distribution $T$ belongs to $\Dlppr (M)$ if
and only if $\alpha * T\in L^p(M)$ for every $\alpha\in\Dcal(M)$.
 \label{th:thconv}
\end{theorem}
{\bf Proof:} Assume first that $1<p\leq\infty$ and let
  $q$ be such that $1/p+1/q=1$. Let $T$ be in $\Dlppr (M)$.
  The set
  $$ B=\{ \phi\in C_c(M): \; ||\phi||_q\leq 1\}$$
 is dense in the unit ball of $L^q(M)$ since
 $1\leq q<\infty$ and $\mu$ is a Radon measure.
 Let $\alpha\in\Dcal(M)$ be fixed. Then by the Young inequality for
 every $\phi\in B$ we get
 $$||\Delta^k(\hat\alpha * \phi)||_q=
   ||\Delta^k\hat\alpha * \phi||_q\leq
   ||\Delta^k\hat\alpha||_1 ||\phi||_q\leq ||\Delta^k\hat\alpha||_1,$$
 where $\hat\alpha=(\widehat{\alpha^\sharp})^\flat$ and
 $\hat{\beta}(g)=\beta(g^{-1})$ for $g\in C_c(G)$.
 This implies that the set
 $\{\hat\alpha*\phi, \phi\in B\}$ is bounded in $\Dlq(M)$ and,
 by using Lemma \ref{l:facts}, we get that
 $$ \langle \alpha *T,\phi\rangle=
    \langle \alpha^\sharp *T^\sharp,\phi^\sharp\rangle=
    \langle T^\sharp,\widehat{\alpha^\sharp} *\phi^\sharp\rangle=
    \langle T,(\widehat{\alpha^\sharp})^\flat *\phi\rangle=
    \langle T, \hat\alpha *\phi \rangle$$
 is bounded if we let $\phi\in B$. It follows that
 $\sup_{\phi\in B} |\langle \alpha * T,\phi\rangle|$ are bounded
 and hence $\alpha*T$ extends to a linear continuous functional on $L^q(M)$,
 which means that $\alpha*T\in L^p(M)$.

 Conversely, suppose that $\alpha*T\in L^p(M)$ for all $\alpha\in\Dcal(M)$.
 For a fixed $\alpha\in\Dcal(M)$,
 $$ \langle \hat\phi *T,\hat\alpha\rangle=\langle\alpha *T,\phi\rangle$$
  are bounded for $\phi\in B$ and, therefore,
 $\hat\phi *T$ are bounded in $\Dcalpr(M)$ for $\phi\in B$.
 Let us denote by $k(m)$ the smallest number $k$ for which
 the fundamental solution $K$ for
 $\Delta^k$ is in $C^m(M)$. Let $F$ be a compact neighborhood of
 the origin of $M$ and let $\gamma\in\Dcal_F(M)$ be equal to one
 in a neighborhood of the origin and supported in $F$. Then
 formula (\ref{eq:fs2}) with $P=\Delta^k$ implies that
  \begin{equation}
   T=\delta * T=\Delta^k(\gamma K *T)-\beta * T,
   \label{eq:eqdec}
  \end{equation}
 where $\beta *T\in L^p(M)$ according to our assumption and
 $\gamma K \in C^m(M)$ is supported in $F$.

 Now we will use the property that bounded sets of distributions
 are equicontinuous on compact sets in $\Rn$ (see \cite{Thomas}), which means
 that if $D^\prime\subset \Dcalpr(\Rn)$ is bounded and $N\subset\Rn$
 is compact, then there exists an integer $m$ such that
  $$|\langle T,\psi\rangle|\leq C \max_{|\nu|\leq m}
     || D^\nu\psi ||_q,\;\forall \psi\in\Dcal_N,\;
     \forall T\in D^\prime.$$
 This estimate is a direct consequence of the uniform boundedness
 principle in the Fr\'{e}chet space $\Dlq$.
 Because $F$ is compact, it follows that distributions
 $\hat\phi *T$ are equicontinuous
 on $F$, since we can apply Theorem \ref{th:th1} to their localizations
 in $\Rn$. Therefore $\hat\phi *T$ extend to $\Dcal^{(m)}_F(M)$
 for all $\phi\in B$ and the extensions satisfy
 $$\sup_{\phi\in B} |\langle \alpha * T,\phi\rangle|<\infty,$$
 in particular if we take $\alpha=\gamma K$.
  It follows that $\gamma K * T\in L^p(M)$ since
  $\gamma K\in \Dcal^{(m)}_F(M)$ if $K$ is a
  fundamental solution for $\Delta^{k(m)}$.
  The statement follows now from Theorem \ref{th:threpr1}
  and decomposition (\ref{eq:eqdec}).

 The case of $p=1$ is similar, implying that $\gamma K * T$ must be
 a Radon measure. A standard additional closed graph
 argument shows that
 it is actually in $L^\infty (M)$.
 Finally we generalize some properties of convolutions.

 \begin{theorem} The following holds:
   \begin{itemize}
     \item[{\rm (i)}] Let $T\in \Dlppr(M)$ and $\phi\in \Dlq(M)$.
      Then $\phi T\in\Dlrpr(M)$ if $r\geq 1$ and $1/r\leq 1/p+1/q$.
     \item[{\rm (ii)}] Let $T\in \Dlppr(M)$ and
           $S\in \Dlqpr(M)$, where $1/p+1/q\geq 1$.
      Then $T*S\in\Dlrpr(M)$ if $1/r=1/p+1/q-1$.
     \item[{\rm (iii)}] Let $T\in\Dlppr(M)$ and
       $\phi\in\Dlq(M)$, where $1/p+1/q\geq 1$.
       Then $\phi*T\in \Dlr(M)$ if $1/r=1/p+1/q-1$.
   \end{itemize}
   \label{th:thests}
 \end{theorem}
 Note that the mappings induced by (i) and (iii), are separately
 continuous in the corresponding spaces, and the mapping in (ii)
 is continuous. The proof of Theorem \ref{th:thests} consists
 of application of Theorem \ref{th:threpr1} to the
 spaces in Theorem \ref{th:thests}.
 The statement follows from the corresponding
 properties of the convolution (\cite{Hew}). We omit the details
 since they are similar to those in Section VI in \cite{Schwartz}
 after we reduce the problems to the Lebesgue spaces by Theorem
 \ref{th:threpr1}.


\medskip

 \noindent
{\sc \small Mathematics Department, Imperial College \\
 180 Queen's Gate, London SW7 2AZ, UK} \\
{\em E-mail address:} ruzh@ic.ac.uk

\medskip
\noindent
\sc
Eingegangen am 5.November 2001

  \end{document}